\documentclass[reqno,11pt]{amsart}

\usepackage{amsmath, amsthm, amssymb, amsfonts, amsthm,tikz}
\usepackage{amsxtra, amscd, geometry, graphicx,epic,eepic}
\usepackage[all,cmtip]{xy}
\usepackage{mathrsfs}
\usepackage{hyperref}
\usepackage{graphicx, nicefrac}
\usepackage{enumerate}
\usepackage{bbm}
\usepackage{holtpolt}

\makeatletter
\newsavebox{\@brx}
\newcommand{\llangle}[1][]{\savebox{\@brx}{\(\m@th{#1\langle}\)}%
  \mathopen{\copy\@brx\kern-0.5\wd\@brx\usebox{\@brx}}}
\newcommand{\rrangle}[1][]{\savebox{\@brx}{\(\m@th{#1\rangle}\)}%
  \mathclose{\copy\@brx\kern-0.5\wd\@brx\usebox{\@brx}}}
\makeatother

\theoremstyle{definition}
\newtheorem*{Ex1*}{Example 1}
\newtheorem*{Ex2*}{Example 2}
\newtheorem*{Ex3*}{Example 3}
\newtheorem*{Ex4*}{Example 4}
\newtheorem*{Cor1*}{Corollary 1}
\newtheorem*{Cor2*}{Corollary 2}
\newtheorem*{Cor3*}{Corollary 3}
\newtheorem*{Cor4*}{Corollary 4}
\newtheorem*{Prop1*}{Proposition 1}
\newtheorem*{Prop2*}{Proposition 2}
\newtheorem*{Lem1*}{Lemma 1}
\newtheorem*{Lem2*}{Lemma 2}
\newtheorem*{Lem3*}{Lemma 3}
\newtheorem*{Thm1*}{Theorem 1}

\newtheorem*{rmk*}{Remark}

\newcommand{\RNum}[1]{\uppercase\expandafter{\romannumeral #1\relax}}

\newcommand{\N}{\mathbb{N}}
\newcommand{\R}{\mathbb{R}}
\newcommand{\Z}{\mathbb Z}
\newcommand{\C}{\mathbb C}
\newcommand{\Q}{\mathbb Q}
\newcommand{\GG}{\mathcal G}

\title[On the Gauss-Kuzmin-L\' evy problem for nearest integer continued fractions]{On the Gauss-Kuzmin-L\' evy problem for nearest integer continued fractions}
\author{Florin P. Boca and Maria Siskaki}

\date{March 31, 2024}

\address{FPB: Department of Mathematics, University of Illinois at Urbana-Champaign, Urbana, IL 61801}
\address{E-mail: fboca@illinois.edu}

\address{MS: Department of Mathematics, University of Illinois at Urbana-Champaign, Urbana, IL 61801}

\curraddr{Department of Mathematics, Yale University, New Haven, CT 06511}

\address{E-mail:  maria.siskaki@yale.edu}

\begin{document}

\begin{abstract}
This note provides an effective bound in the Gauss-Kuzmin-L\' evy problem for some Gauss type shifts associated with
nearest integer continued fractions, acting on the interval $I_0=[0,\frac{1}{2}]$ or $I_0=[-\frac{1}{2},\frac{1}{2}]$.
We prove asymptotic formulas $\lambda (T^{-n}I) =\mu(I)(\lambda ( I_0) +O(q^n))$
for such transformations $T$, where $\lambda$ is the Lebesgue measure on $\mathbb R$,
$\mu$ the normalized $T$-invariant Lebesgue absolutely continuous measure,
$I$ subinterval in $I_0$, and $q=0.288$ is smaller than the Wirsing constant $q_W=0.3036\ldots$
\end{abstract}

\maketitle

\section{Introduction}
The regular continued fraction establishes a one-to-one correspondence between the set  of
infinite words with letters in the alphabet $\N$ and the set $[0,1]\setminus {\mathbb Q}$:
$$
(a_1,a_2,a_3,\ldots) \mapsto \polter{1}{a_1} +\polter{1}{a_2} +\polter{1}{a_3} +\cdots
$$

The Gauss shift $\GG$ acts on $\N^\N$ by $\GG (a_1,a_2,a_3,\ldots )=(a_2,a_3,a_4,\ldots)$, and on $[0,1]$ by
$\GG (x)=\{ \frac{1}{x}\} =\frac{1}{x}-\lfloor \frac{1}{x}\rfloor$ if $x\neq 0$ and $\GG (0)=0$.
Gauss discovered that the probability measure $d\mu =\frac{dx}{(1+x)\log 2}$ is $\GG$-invariant on $[0,1]$ and stated in his diary (October 25, 1800) that
\begin{equation}\label{eq1}
\lambda (\GG^{-n} [0,x]) \sim \mu ([0,x])=\frac{\log(1+x)}{\log 2} ,\quad \forall x\in [0,1] \text{ as $n\rightarrow \infty$,}
\end{equation}
where $\lambda$ denotes the Lebesgue measure on $\R$.
In a 1812 letter to Laplace (\cite{Ga}, see also Appendix III of \cite{Usp}), Gauss raised  the problem of providing an effective version of \eqref{eq1} and estimate the error
$$
E_n (x)=\lambda (\GG^{-n}[0,x]) -\mu ([0,x]),\quad x\in [0,1],\ n\rightarrow \infty .
$$
The problem was thoroughly investigated  much later, with significant contributions by
Kuzmin \cite{Kuz} and L\' evy \cite{Lev}. Kuzmin proved that, uniformly in $x$, $E_n(x)=O(q^{\sqrt{n}})$ for some $q\in (0,1)$, while L\' evy proved that
$E_n(x)=O(q^n)$ with $q< 0.7$. The breakthrough result of Wirsing \cite{Wir} proved that
\begin{equation}\label{eq2}
E_n(x) =\psi (x) q_W^n +O(q_1^n),
\end{equation}
with $q_W =0.303663\ldots$ denoting the Wirsing (optimal) constant, $0<q_1 <q_W$ and $\psi$ some real analytic function on $[0,1]$.
The spectral approach due to Babenko \cite{Bab} and Mayer-Roepstorff \cite{MR} provided a complete solution to the
problem, showing that the restriction of the Perron-Frobenius operator of $\GG$ to some Hardy space on the right half-plane
$\operatorname{Re} z > -\frac{1}{2}$ is similar to a self-adjoint trace-class operator with explicit kernel,
and thus the expression of $E_n(x)$ in \eqref{eq2} can be completed to the eigenfunction expansion of this compact operator.
A detailed discussion of the Gauss problem with complete proofs can be found in the monograph \cite{IK}.

It is natural to study the analogue of the Gauss problem for other classes of continued fractions.
This note takes an elementary look at the situation of the nearest integer continued fraction (\emph{NICF}), originally considered in Minnigerode's work
on the Pell equation \cite{Min} and furthered by Hurwitz \cite{Hur}.
NICF provides a better rate of approximation than the regular continued fraction. Actually, each nearest integer convergent of an irrational number is a
regular continued fraction convergent of that number \cite{Ada,Wil}. Other Diophantine approximation properties, such as analogues of Vahlen's theorem,
were studied in \cite{JK,Ton}. Analogues of the Gauss problem for other types of continued fractions have been recently studied in \cite{LS1,LS2,Sun}.

In this paper we denote $G=\frac{\sqrt{5}+1}{2}$, $g=\frac{\sqrt{5}-1}{2}$, and employ the equalities $G-1=g$, $G+1=G^2$ and
$(2-G)(G+1)=1$.

The NICF can appear in various guises. We will consider three possible situations, as follows:

(A) The folded NICF map $T:[0,\frac{1}{2}]\longrightarrow [0,\frac{1}{2}]$ defined by $T(0)=0$, and, for $x \neq 0$, by
\begin{equation}
T(x)=\bigg|\,\frac{1}{x}-\bigg\lfloor \frac{1}{x}+\frac{1}{2}\bigg\rfloor \,\bigg|
=\begin{cases} \lvert \frac{1}{x}-k\rvert & \mbox{\rm if $\frac{2}{2k+1} \leq x \leq \frac{2}{2k-1},\ k\geq 3$} \\
\frac{1}{x} -2 & \mbox{\rm if $\frac{2}{5} \leq x \leq \frac{1}{2}$}
\end{cases} ,
\end{equation}
is continuous on $(0,\frac{1}{2}]$. For every $x\in [0,\frac{1}{2}] \setminus \Q$, let $a_1=a_1(x):=\lfloor x+\frac{1}{2}\rfloor \geq 2$,
$e_1=e_1 (x):=\operatorname{sign} (\frac{1}{x}-a_1)\in \{ \pm 1\}$. They satisfy $a_1+e_1 \geq 2$.
Note that $T(x)=e_1 (\frac{1}{x}-a_1)=\lvert \frac{1}{x}-a_1\rvert$, $\forall x\in [0,\frac{1}{2}]\setminus \Q$.
Taking $a_i=a_i(x):=a_1 (T^{i-1}(x))$, $e_i=e_i(x):=e_1(T^{i-1}(x))$ if $i\geq 2$, every irrational $x\in [0,\frac{1}{2}]$ is represented as
$$
x= \polter{1}{a_1} +\polter{e_1}{a_2}+\polter{e_2}{a_3}+\cdots =: [(a_1,e_1),(a_2,e_2),(a_3,e_3),\ldots] ,
$$
with $a_i\geq 2$, $e_i\in \{ \pm 1\}$ and $a_i+e_i \geq 2$.

The map $T$ is called a Gauss type shift because it acts on $[0,\frac{1}{2}]\setminus \Q$ by shifting the digits $(a_i,e_i)$:
$$T
([(a_1,e_1),(a_2,e_2),\ldots ])=[(a_2,e_2),(a_3,e_3),\ldots ].
$$

According to Lemma 1 below, the probability measure
$$
d\mu =\frac{1}{\log G} \bigg( \frac{1}{G+x}+\frac{1}{G+1-x}\bigg) dx
$$
is $T$-invariant.

\medskip

(B) The odd map $T_o :[-\frac{1}{2},\frac{1}{2}]\longrightarrow [-\frac{1}{2},\frac{1}{2}]$, investigated by Nakada, Ito and Tanaka \cite{NIT} and defined by $T_o(0)=0$, and, for $x \neq 0$, by
\begin{equation}
T_o (x)=\frac{1}{x}-\bigg\lfloor \frac{1}{x}+\frac{1}{2}\bigg\rfloor=\begin{cases}
\frac{1}{x}-k\operatorname{sgn} x & \text{if $\frac{2}{2k+1} < \lvert x\rvert \leq \frac{2}{2k-1}$, $k\geq 3$} \\
\frac{1}{x}-2\operatorname{sgn} x & \text{if $\frac{2}{5} < \lvert x\rvert\leq \frac{1}{2}$}
\end{cases} =-T_o(-x),
\end{equation}
represents the Gauss shift associated with the continued fraction expansion
$$x=\polter{1}{b_1}+\polter{1}{b_2}+\polter{1}{b_3}+\cdots =: [b_1,b_2,b_3,\ldots]
$$
of irrationals in $[-\frac{1}{2},\frac{1}{2}]$, with digits $b_i=b_i(x) \in \Z$ given by
$b_1 :=\lfloor \frac{1}{x}+\frac{1}{2}\rfloor$, $b_i:=b_1(T_o^{i-1}(x)) $ if $i\geq 2$. Then
$\lvert b_i \rvert \geq 2$, $b_i=2 \Longrightarrow b_{i+1}\geq 2$, and
$b_i=-2 \Longrightarrow b_{i+1}\leq -2$.
Indeed, it is plain that $T_o (x)=\frac{1}{x}-b_1$, $\forall x\in [-\frac{1}{2},\frac{1}{2}]\setminus {\mathbb Q}$, and
$T_o ([b_1,b_2,\ldots ])=[b_2,b_3,\ldots ]$.
As shown in \cite{NIT}, the probability measure
$$
d\mu_o =\frac{1}{2\log G} \bigg( \frac{1}{G+\lvert x\rvert} +\frac{1}{G+1-\lvert x\rvert} \bigg)dx
$$
is $T_o$-invariant.

The identity
$
-\polter{1}{b_1} -\polter{1}{b_2} -\polter{1}{b_3} -\cdots = \polter{1}{-b_1}+\polter{1}{b_2}+\polter{1}{-b_3}+\polter{1}{b_4}+\cdots
$
shows that $T_o$ can also be viewed as the Gauss shift generated by the NICF expansion
$b_0 -\polter{1}{b_1}-\polter{1}{b_2} -\cdots$ considered in \cite{Ada,Hur}.

\medskip

(C) The even map $T_e:[-\frac{1}{2},\frac{1}{2}]\longrightarrow [-\frac{1}{2},\frac{1}{2}]$, considered by Rieger \cite{Ri1,Rie} and defined by $T_e(0)=0$, and, for $x \neq 0$, by
$$
T_e (x):=\frac{1}{\lvert x\rvert} -\bigg\lfloor \frac{1}{\lvert x\rvert} +\frac{1}{2}\bigg\rfloor=
\begin{cases} \frac{1}{\lvert x\rvert} -k & \text{if $\frac{2}{2k+1} < \lvert x\rvert \leq \frac{2}{2k-1}$, $k\geq 3$} \\
\frac{1}{\lvert x\rvert}-2 & \text{ if $\frac{2}{5} < \lvert x\rvert \leq \frac{1}{2}$}
\end{cases} =T_e (-x),
$$
generates the NICF expansion
$$
x=\polter{e_1}{a_1} +\polter{e_2}{a_2}+\polter{e_3}{a_3}+\cdots =: [\![ (a_1,e_1),(a_2,e_2),(a_3,e_3),\ldots ]\!]
$$
of irrationals in $[-\frac{1}{2},\frac{1}{2}]$,
with digits $a_1=a_1(x):=\lfloor \frac{1}{\lvert x\rvert}+\frac{1}{2}\rfloor$, $e_1=e_1(x):=\operatorname{sign}  (\frac{1}{\lvert x\rvert}-a_1)$,
$a_i:=a_1(T_e^{i-1}(x))$, $e_i :=e_1 (T_e^{i-1}(x))$ if $i\geq 2$ satisfying
$a_i\geq 2$, $e_i \in \{ \pm 1\}$, $a_i+e_{i+1} \geq 2$.
This NICF expansion is also considered in \cite{JK,Per,Ton}.

We have $T_e (x) =\frac{e_1}{x}-a_1=\frac{1}{\lvert x\rvert}-a_1$ and
$T_e ([\![ (a_1,e_1),(a_2,e_2),\ldots ]\!])=[\![(a_2,e_2),(a_3,e_3),\ldots ]\!]$,
so $T_e$ is the Gauss shift associated with this NICF expansion.

The map $T_e$ coincides with Nakada's map $f_{1/2}$ \cite{Nak}. In particular, the probability measure
$$
d\mu_e =h_e(x) dx ,\quad  h_e(x)=\frac{1}{\log G} \begin{cases} \frac{1}{G+x} & \text{ if $0<x<\frac{1}{2}$} \\
\frac{1}{G+x+1} & \text{ if $-\frac{1}{2} <x<0$}
\end{cases}
$$
is $T_e$-invariant.

\begin{figure}
\includegraphics[scale=0.5]{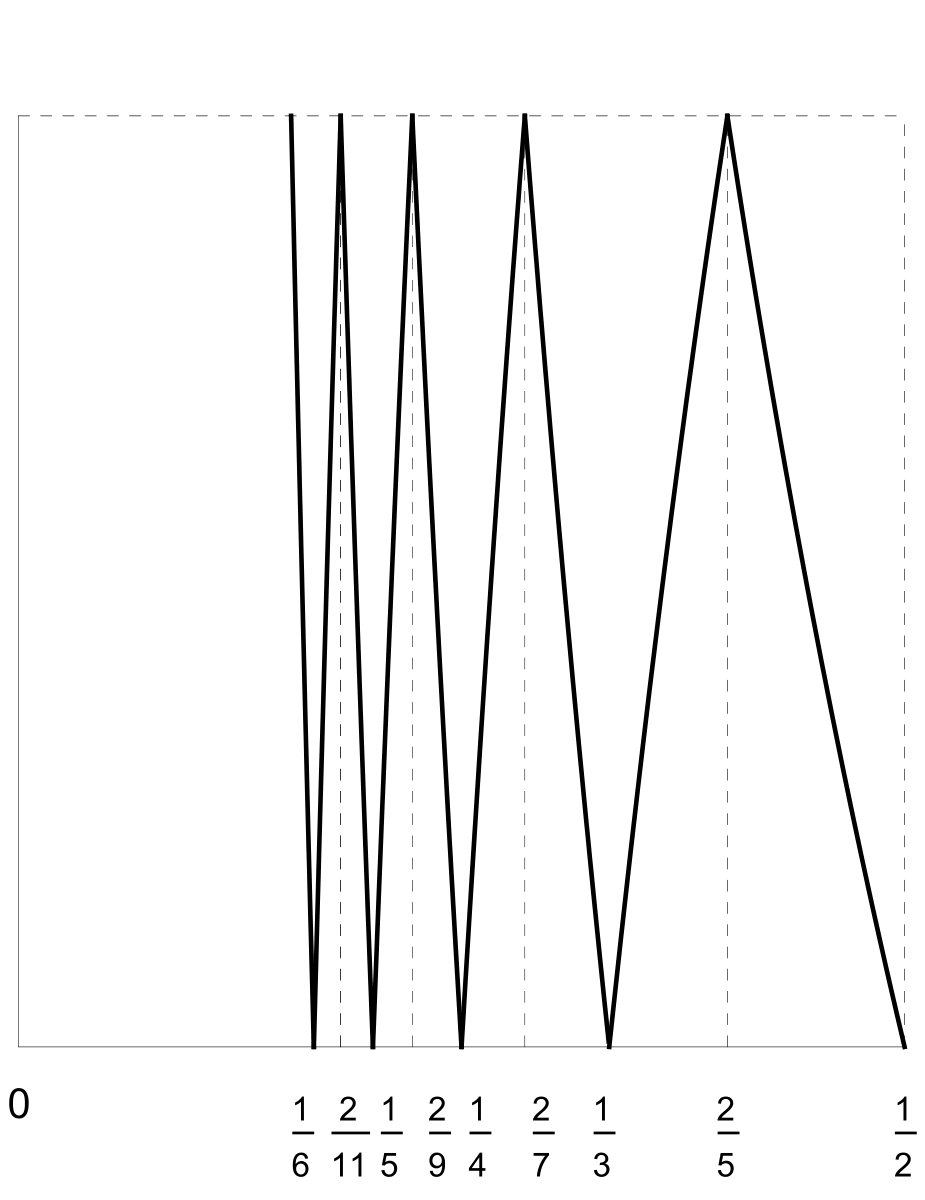}\qquad
\includegraphics[scale=0.52]{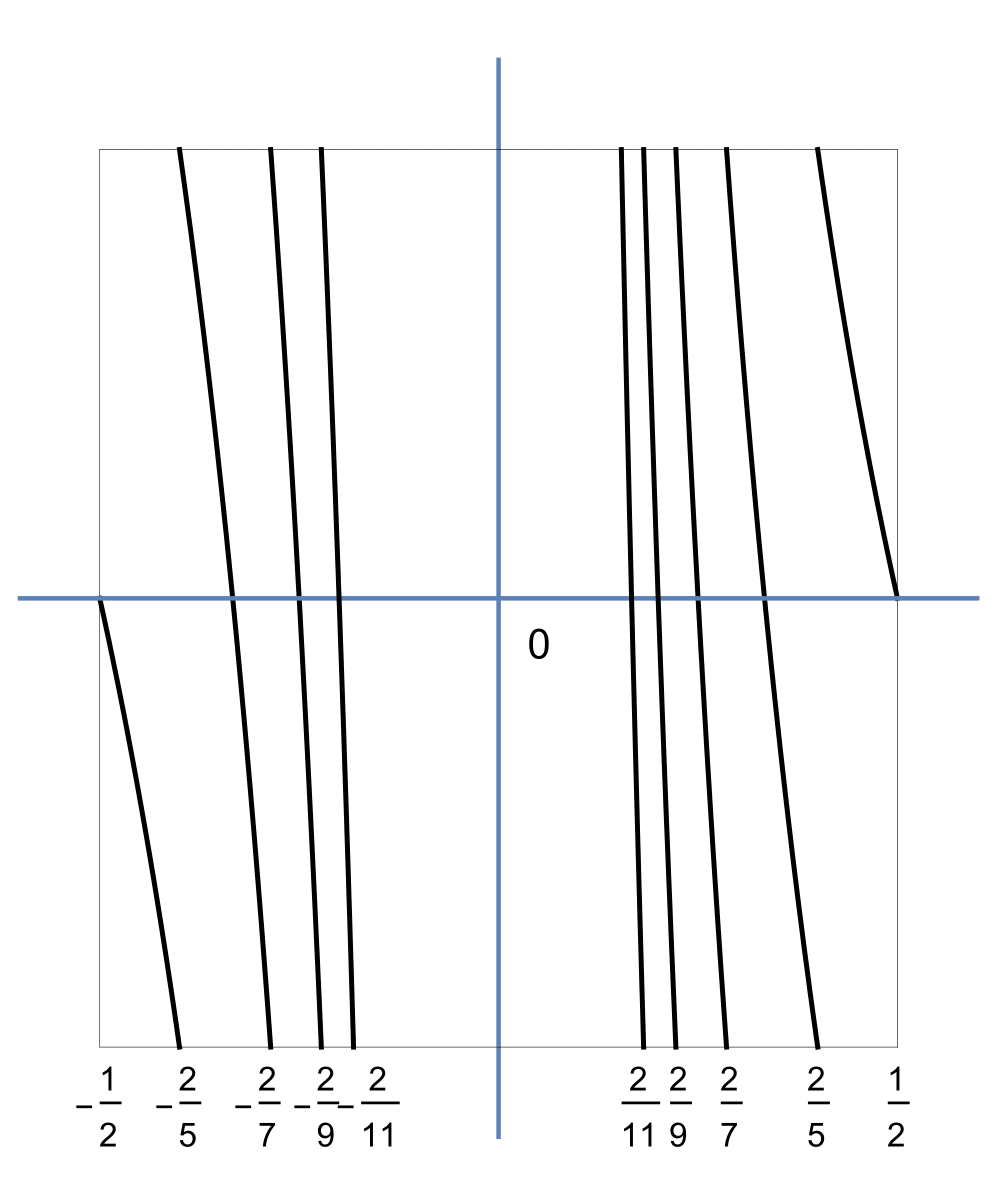}
\caption{Graphs of the maps $T$ and $T_o$.}
\label{Figure1}
\end{figure}

The main result of this note provides quantitative estimates for the analogue of the Gauss-Kuzmin-L\' evy problem in the
situations of the Gauss type shifts $T$, $T_o$ and $T_e$, as follows:

\begin{Thm1*}
\begin{itemize}
\item[(i)]
With $q=0.288$, for every Borel set $E\subseteq [0,\frac{1}{2}]$,
$$
\lambda (T^{-n} E)=\frac{1}{2}\,\mu(E) +O(\mu(E) q^n).
$$
\item[(ii)]
With $q=0.288$, for every Borel set $E\subseteq [-\frac{1}{2},\frac{1}{2}]$,
$$
\lambda (T_o^{-n} E)=\mu_o (E)+O (\mu_o (E) q^n).
$$
\item[(iii)]
With $q=0.234$, for every Borel set $E\subseteq [-\frac{1}{2},\frac{1}{2}]$,
$$
\lambda (T_e^{-n} E)=\mu_e (E) +O(\mu_e (E) q^n).
$$
\end{itemize}
\end{Thm1*}

The estimate in (ii) improves upon $q=g^2 \approx 0.382$ obtained in \cite[Thm.2.1(ii)]{NIT}.
The estimate in (iii) improves upon $q=\frac{2}{3}$ obtained in \cite{Ri1}. Note that
$q=0.288$ is smaller that the Wirsing constant $q_W =0.3036\ldots$.

To prove Theorem 1, we perform an elementary analysis of the Perron-Frobenius operators associated
to the transformations $T$ and $T_e$ with respect to their invariant Lebesgue absolutely continuous measures along the line of \cite{Ri1}.

In \cite{IKa} and \cite{Pop}, the authors investigated a problem similar to (iii). However, their
transition operator $U$ coincides with the Perron-Frobenius operator associated to the \emph{dual} of the NICF Gauss map, rather than the NICF Gauss map itself.
This dual is the folded Hurwitz transformation
$S$, which acts on $[0,g]$ by $S(0)=0$ and
$$
S(x)=\bigg| \frac{1}{x}-i\bigg| \quad \text{if $\displaystyle \frac{1}{i+g}<x\leq \frac{1}{i+g-1}$},\ i\geq 2,
$$
with $S$-invariant probability measure
$$
d \nu =k(x) dx,\quad k(x)=\frac{1}{\log G}\begin{cases} \frac{1}{2+x}+\frac{1}{2-x} & \text{if $x\in [0,g^2)$} \\
\frac{1}{2+x} & \text{if $x\in [g^2,g]$}.
\end{cases}
$$

We also provide some estimates on the rate of mixing of the map $T$.

\begin{Cor1*}
With $q=0.288$, for any Borel set $E\subseteq [0,\frac{1}{2}]$ and any $T$-cylinder $F$,
$$
\mu (T^{-n} E \cap F)=\mu (E)\mu(F)+O_F (q^n).
$$
\end{Cor1*}

\begin{Cor2*}
With $q=0.288$, for any Borel set $E\subseteq [-\frac{1}{2},\frac{1}{2}]$ \emph{symmetric} with respect to the origin, and any $T_o$-cylinder $F$,
$$
\mu_o (T_o^{-n} E \cap F) =\mu_o(E)\mu_o(F)+O_F (q^n).
$$
\end{Cor2*}

Corollary 2 was proved with $q=g^2$ in \cite[Thm.2.1(iii)]{NIT} without assuming $E$ symmetric.

An analogue of Corollaries 1 and 2 will be discussed at the end of Section 3.

%\begin{Cor3*}
%With $q=0.234$, for any Borel set $E\subseteq [-\frac{1}{2},\frac{1}{2}]$ and any $T_e$-cylinder $F$,
%$$
%\mu_e(T_e^{-n} E \cap F)=\mu_e(E) \mu_e(F)+O_F (q^n).
%$$
%\end{Cor3*}

\section{The folded NICF map $T$ and the Nakada-Ito-Tanaka map $T_o$}
The folded NICF can be obtained as a particular example of a folded Japanese continued fraction,
investigated by Moussa, Cassa and Marmi \cite{MCM}.
The following lemma follows by taking $\alpha=\frac{1}{2}$ in \cite[Thm.15]{MCM}, or it can be verified directly through a plain calculation.

\begin{Lem1*}
The probability measure $d\mu =C h(x) dx$, with
$$
h(x):=\frac{1}{G+x}+\frac{1}{G+1-x} ,\qquad C=\frac{1}{\log G},
$$
is $T$-invariant.
\end{Lem1*}

Denote
$$
W  :=\{ (k,1) : k\geq 2\} \cup \{ (k,-1): k\geq 3\} , \qquad
w_{(k,e)}(y)  :=\frac{1}{k+ey} .
$$
Following \cite[Sect.2.3]{KMS},
the Perron-Frobenius (Ruelle) operator $P =\widehat{T}_\lambda$ of $T$ with respect to the Lebesgue measure $\lambda$ acts on $L^1([0,\frac{1}{2}],\lambda)$  by
\begin{equation}\label{eq_P}
\begin{aligned}
(P f)(y) & =\sum\limits_{x\in T^{-1} y} f(x(y))\lvert x^\prime (y)\rvert
=\sum\limits_{(k,e)\in W} w_{(k,e)}^2 (y) f(w_{(k,e)} (y)) \\ & =
\sum\limits_{k\geq 2} w_{(k,1)}^2 (y) f(w_{(k,1)}(y)) +\sum\limits_{k\geq 3} w_{(k,-1)}^2 (y) f(w_{(k,-1)} (y)) .
\end{aligned}
\end{equation}
The Perron-Frobenius (transfer) operator $U:=\widehat{T}_\mu$ of $T$
with respect to the invariant measure $\mu$ acts on $L^1([0,\frac{1}{2}],\mu)$ by
\begin{equation}\label{eq1.1}
U=M_H PM_H^{-1} =M_H P M_h ,
\end{equation}
where $M_H$ denotes the operator of multiplication by $H:=\frac{1}{h}$. Since $\mu$ is a $T$-invariant measure, one has  $U1=1$.
One can also consider $U$ as the transpose (dual) of the Koopman
operator defined by $K_T f:=f\circ T$.

The equalities \eqref{eq_P}, \eqref{eq1.1} and
$$
\begin{aligned}
w_{(k,e)}^2 (y) h(w_{(k,e)}(y)) & =\frac{1}{(k+ey)^2} \bigg( \frac{1}{G+\frac{1}{k+ey}} +\frac{1}{G+1-\frac{1}{k+ey}}\bigg) \\
& =  \frac{1}{k+G-2+ey} -\frac{1}{k+G-1+ey}
\end{aligned}
$$
lead to
\begin{equation}\label{eq1.2}
(Uf)(y) = \sum\limits_{(k,e)\in W} P_{(k,e)} (y) f(w_{(k,e)}(y)) ,
\end{equation}
with
\begin{equation}\label{eq P_(k,e)}
 P_{(k,e)}(y)=H(y) \bigg( \frac{1}{k+G-2+ey} -\frac{1}{k+G-1+ey}\bigg) \geq 0.
\end{equation}

The equalities \eqref{eq1.2} and \eqref{eq P_(k,e)} show that $U(C[0,\frac{1}{2}]) \subseteq C[0,\frac{1}{2}]$ and
$U(C^1[0,\frac{1}{2}]) \subseteq C^1[0,\frac{1}{2}]$.

Since $U1=1$, the weights
$
P_{(k,e)}
$
satisfy
\begin{equation}\label{eq1.3}
\sum\limits_{(k,e)\in W} P_{(k,e)} (y)=1 \quad \text{and}\quad \sum\limits_{(k,e)\in W} P^\prime_{(k,e)} (y)=0,\quad \forall y\in [0,\tfrac{1}{2}] .
\end{equation}
The first identity in \eqref{eq1.3} allows us to write
\begin{equation*}
(Uf)(y)=f(\tfrac{1}{4})+\sum\limits_{(k,e)\in W} P_{(k,e)}(y) (f(w_{(k,e)}(y) -f(\tfrac{1}{4})).
\end{equation*}

Inserting $A=k+ey$, $B=G-2=-\frac{1}{G^2}$, which satisfy $\frac{1}{B}-\frac{1}{B+1}=-G^3$ and $\frac{1}{B}+\frac{1}{B+1}=-1$, in the identity
$$
\begin{aligned}
\frac{1}{A^2} \bigg( \frac{1}{A+B} -\frac{1}{A+B+1}\bigg) & =\bigg( \frac{1}{B}-\frac{1}{B+1}\bigg) \frac{1}{A^2}
-\bigg( \frac{1}{B^2} -\frac{1}{(B+1)^2}\bigg) \frac{1}{A} \\ & \qquad
+\frac{1}{B^2}\cdot \frac{1}{A+B} -\frac{1}{(B+1)^2} \cdot \frac{1}{A+B+1} ,
\end{aligned}
$$
we infer
\begin{equation}\label{eq1.4}
\begin{aligned}
\frac{P_{(k,e)}(y)}{(k+ey)^2} & =H(y) \bigg( -\frac{G^3}{(k+ey)^2} -\frac{G^3}{k+ey} +\frac{G^3+G^2}{k+G-2+ey} -\frac{G^2}{k+G-1+ey}\bigg) \\
& = H(y) \bigg( -\frac{G^3}{(k+ey)^2} +\frac{G}{(k+ey)(k+G-2+ey)} \\ & \qquad \qquad \qquad  +G^2 \Big( \frac{1}{k+G-2+ey}-\frac{1}{k+G-1+ey}\Big) \bigg).
\end{aligned}
\end{equation}

\begin{Prop1*}
$
\| (Uf)^\prime \|_\infty \leq 0.288\| f^\prime\|_\infty$, $\forall f\in C^1 [0,\tfrac{1}{2}] .
$
\end{Prop1*}

\proof[Proof]
Employing \eqref{eq1.2}, \eqref{eq1.3}, the Mean Value Theorem, and $\lvert w_{(k,e)} (y)-\frac{1}{4}\rvert \leq \frac{1}{4}$,  we can write
\begin{equation}
\begin{aligned}
\lvert (Uf)^\prime (y)\rvert & \leq \sum\limits_{(k,e)\in W} \frac{P_{(k,e)}(y)}{(k+ey)^2} \lvert f^\prime (w_{(k,e)} (y))\rvert
+\sum\limits_{(k,e)\in W} \lvert P^\prime_{(k,e)}(y) (f(w_{(k,e)} (y))-f(\tfrac{1}{4}))\rvert \\
& \leq (S_I(y)+S_{II}(y)) \| f^\prime\|_\infty ,\quad \forall y\in [ 0,\tfrac{1}{2}],
\end{aligned}
\end{equation}
with
$$
S_I(y)=\sum\limits_{(k,e)\in W} \frac{P_{(k,e)}(y)}{(k+ey)^2},\qquad S_{II}(y)=\frac{1}{4} \sum\limits_{(k,e)\in W} \lvert P^\prime_{(k,e)}(y)\rvert .
$$

The identity
$$
\sum\limits_{k\in\Z} \frac{1}{(k+z)^2}=\frac{\pi^2}{\sin^2 (\pi z)},\quad \forall z\in \C \setminus \Z ,
$$
leads to
$$
\Phi_1(y):=\sum\limits_{(k,e)\in W} \frac{1}{(k+ey)^2} =\frac{\pi^2}{\sin^2 (\pi y)} -\frac{1}{y^2}-\frac{1}{(1+y)^2}-\frac{1}{(1-y)^2}-\frac{1}{(2-y)^2} ,
\quad \forall y\in (0,\tfrac{1}{2}] .
$$
Note that $\Phi_1 (0^+)=\frac{\pi^2}{3}-\frac{9}{4}$.

We also have
$$
\begin{aligned}
& H(y) \sum\limits_{(k,e)\in W} \bigg( \frac{1}{k+G-2+ey}-\frac{1}{k+G-1+ey}\bigg) \\ &
 =H(y) \sum\limits_{k\geq 2} \bigg( \frac{1}{k+G-2+y}-\frac{1}{k+G-1+y}\bigg)
+H(y) \sum\limits_{k\geq 3} \bigg( \frac{1}{k+G-2-y}-\frac{1}{k+G-1-y}\bigg) \\
& =H(y) \bigg( \frac{1}{G+y}+\frac{1}{G+1-y}\bigg) =H(y) h(y)=1.
\end{aligned}
$$
Combining \eqref{eq1.4} with the last two equations above and using $H(y)=\frac{(G+y)(G+1-y)}{G^3}$ we find
$$
S_I(y) =-(G+y)(G+1-y) \Phi_1 (y) +GH(y)\Phi_2(y)+G^2,
$$
where
$$\Phi_2(y) :=\sum\limits_{(k,e)\in W} \frac{1}{(k+ey)(k+G-2+ey)} .
$$

Numerically, Mathematica gives
$$
\begin{aligned}
S_I(y) \leq S_I(0^+) =G^2-G^3\bigg( \frac{\pi^2}{3}-\frac{9}{4}\bigg) +G\Phi_2(0) < 0.097,\quad \forall y\in [0,\tfrac{1}{2}].
\end{aligned}
$$

To bound $S_{II}(y)$, we write $P_{(k,e)} =L_{(k,e)}-L_{(k+1,e)}$, with
$$
\begin{aligned}
L_{(k,e)}(y)& =\frac{H(y)}{k+ey+G-2}=\frac{1}{G^3}\cdot \frac{G^3+y-y^2}{k+ey+G-2}, \\
L^\prime_{(k,e)} (y) & = \frac{1}{G^3} \bigg( \frac{1-2y}{k+ey+G-2} -\frac{e(G^3+y-y^2)}{(k+ey+G-2)^2}\bigg) ,
\end{aligned}
$$
and compute
$$
\begin{aligned}
P^\prime_{(k,e)}(y)& =\frac{1-2y}{G^3} \bigg( \frac{1}{k+ey+G-2} -\frac{1}{k+ey+G-1}\bigg) \\ & \qquad
-\frac{e(G^3+y-y^2)}{G^3} \bigg( \frac{1}{(k+ey+G-2)^2} -\frac{1}{(k+ey+G-1)^2}\bigg).
\end{aligned}
$$
Summing over $(k,e)\in W$, we find
$$
\begin{aligned}
S_{II}(y) & \leq \frac{1-2y}{4G^3} \bigg( \frac{1}{G+y} +\frac{1}{G+1-y}\bigg)
+\frac{(G+y)(G+1-y)}{4G^3} \bigg( \frac{1}{(G+y)^2}+\frac{1}{(G+1-y)^2}\bigg)
\\
& =\frac{1-2y}{4(G+y)(G+1-y)} +\frac{(G+y)(G+1-y)}{4G^3} \bigg( \frac{1}{(G+y)^2}+\frac{1}{(G+1-y)^2}\bigg) \\
& < 0.191,\quad \forall y\in [0,\tfrac{1}{2}].
\end{aligned}
$$
We infer $S_I(y)+S_{II}(y) < 0.097+0.191 =0.288$, $\forall y\in [0,\frac{1}{2}]$.
\qed

\proof[Proof of {\em (i)}and {\em (ii)} in Theorem 1]
(i) Consider $\gamma_n:=U^n H \in C^1 [0,\frac{1}{2}]$. With $h$ as in Lemma 1, $H:=\frac{1}{h}$, and taking $d\nu: =h(x) dx$, we have
\begin{equation}\label{eq1.6}
\lambda (T^{-n} E)=\int_E U^n H\, d\nu =\int_E \gamma_n h\, d\lambda .
\end{equation}
Proposition 1 shows that
$$\| \gamma_n^\prime \|_\infty \leq \| H^\prime\|_\infty q^n <q^n,\quad \forall n\geq 1.
$$
The Mean Value Theorem then yields
$$
\lvert \gamma_n (x)-\gamma_n (0)\rvert \leq q^n x\quad \forall x\in [0,\tfrac{1}{2}],
$$
or equivalently
\begin{equation}\label{eq1.7}
\lvert \gamma_n (x) h(x)-\gamma_n (0) h(x)\rvert \leq q^n xh(x) \leq q^n ,\quad \forall x\in [0,\tfrac{1}{2}], \ \forall n\geq 1.
\end{equation}
Therefore, we have
$$
\begin{aligned}
\frac{1}{2} = & \ \lambda(T^{-n} [0,\tfrac{1}{2}]) =\int_0^{1/2} \gamma_n h\, d\lambda =\gamma_n (0) \int_0^{1/2}h\, d\lambda +O(q^n) \\
& \Longrightarrow\ \gamma_n(0)=\frac{1}{2\int_0^{1/2} h\, d\lambda} +O(q^n) \\
& \stackrel{\eqref{eq1.7}}{\Longrightarrow} \ \gamma_n (x) h(x) =\frac{h(x)}{2\int_0^{1/2} h\, d\lambda} +O(q^n)\\
& \stackrel{\eqref{eq1.6}}{\Longrightarrow}\ \lambda (T^{-n}E) =\frac{\int_E h\, d\lambda}{2\int_0^{1/2} h\, d\lambda} +O(\mu(E) q^n)
=\frac{1}{2}\, \mu(E)+ O(\mu(E) q^n) .
\end{aligned}
$$
In the last line above we also used $\lambda \ll \mu\ll \lambda$.

(ii) The probability measure $\mu_o$ from the introduction is
$T_o$-invariant. Furthermore, the measure $\mu$ is equal to two times the push-forward of  $\mu_o$ under the map
$\lvert\  \rvert :[-\frac{1}{2},\frac{1}{2}]\longrightarrow [0,\frac{1}{2}]$.

Consider a Borel set $E\subseteq [0,\frac{1}{2}]$. We have $T(x)=\lvert T_o (x)\rvert$, $\forall x\in [0,\frac{1}{2}]$, so
$T=\lvert T_o \rvert\, \big\vert_{[0,1/2]}$ and $T^n =\lvert T_o\rvert^n \,\big\vert_{[0,1/2]}$, $\forall n\geq 1$.
Each map $T_o^n$ is odd, so $T_o^{-n} (-E)=-T_o^{-n} E$. This entails
$$
T^{-n} E=\{ x\in [0,\tfrac{1}{2}] : T_o^n (x) \in E \cup (-E)\}
=(T_o^{-n} E \cup (-T_o^{-n} E)) \cap [0,\tfrac{1}{2}] .
$$
The conclusion follows from
$$
\begin{aligned}
\lambda (T^{-n} E) & =\lambda (T_o^{-n} E \cap [0,\tfrac{1}{2}]) +\lambda ((-T_o^{-n} E) \cap [0,\tfrac{1}{2}])
\\ &  =\lambda (T_o^{-n} E \cap [0,\tfrac{1}{2}])+\lambda (T_o^{-n} E \cap [-\tfrac{1}{2},0])
=\lambda (T_o^{-n} E)
\end{aligned}
$$
and Theorem 1 (i), using also $\mu(E)=2\mu_o (E)$.

When $E\subseteq [-\tfrac{1}{2},0]$, we use again $T_o^{-n}(E)=-T_o^{-n}(-E)$ and $\mu_o(-E)=\mu_o(E)$.
\qed

The $T$-cylinders are given by
$\Delta_{[(a_1,e_1)]}:=\{\frac{1}{a_1+e_1 y}:y\in [0,\frac{1}{2}]\}$ and when $r\geq 2$ by
$$
\Delta_{[(a_1,e_1), \ldots ,(a_r,e_r)]} =\bigcap\limits_{i=1}^r T^{-(i-1)} \Delta_{[(a_i,e_i)]}=
\bigg\{\, \polter{1}{a_1} +\polter{e_1}{a_2} + \cdots + \polter{e_{r-1}}{a_r+e_r y} : y\in [ 0,\tfrac{1}{2}]\bigg\} .
$$

\proof[Proof of Corollary 1]
Let $F=\Delta_{[(a_1,e_1), \ldots ,(a_r,e_r)]}$.
We will estimate
\begin{equation}\label{eq1.8}
\mu (T^{-n}E \cap F)=\int_{T^{-n}E} \chi_F \, d\mu =\int_E U^n \chi_F \, d\mu .
\end{equation}
From $\chi_F \circ w_{(k,e)}=\delta_{(k,e),(a_1,e_1)} \chi_{\Delta_{[(a_2,e_2),\ldots,(a_r,e_r)]}}$ and equality \eqref{eq1.2} we infer
$$
U\chi_F =P_{(a_1,e_1)} \cdot \chi_{\Delta_{[(a_2,e_2),\ldots,(a_r,e_r)]}},\quad U^2 \chi_F=(P_{(a_1,e_1)} \circ w_{(a_2,e_2)})\cdot P_{(a_2,e_2)}\cdot
\chi_{\Delta_{[(a_3,e_3),\ldots,(a_r,e_r)]}},
$$
and finally,
\begin{equation}\label{eq1.9}
U^r \chi_F =\prod\limits_{i=1}^r  P_{(a_i,e_i)} \circ w_{(a_{i+1},e_{i+1})} \circ \cdots \circ w_{(a_r,e_r)} =: C_F \in C^1[0,\tfrac{1}{2}],
\end{equation}
where the term corresponding to $i=r$ is just $P_{(a_r,e_r)}(y)$.

Proposition 1 and equality \eqref{eq1.9} entail $\| (U^n \chi_F)^\prime \|_\infty =\| (U^{n-r} C_F)^\prime \|_\infty \ll_F q^n$, $\forall n\geq r$, with $q=0.288$. Applying the Mean Value Theorem  we get
\begin{equation}\label{eq1.10}
\begin{aligned}
\| U^n \chi_F \cdot h -(U^n \chi_F)(0)\cdot h\|_\infty &  =\|U^{n-r} C_F \cdot h-(U^{n-r} C_F)(0) \cdot h\|_\infty \\ &
\leq \|U^{n-r} C_F   -(U^{n-r} C_F)(0) \|_\infty \ll_F q^n,\quad \forall n\geq r.
\end{aligned}
\end{equation}
Integrating on $[0,\frac{1}{2}]$ this yields
$$
\mu (F)=\int_0^{1/2} U^n \chi_F \, d\mu =(U^n \chi_F)(0) +O_F (q^n),\quad \forall n\geq r.
$$
Plugging this back in \eqref{eq1.10} we find
$$
(U^n\chi_F)(x)=\mu(F) +O_F (q^n),\quad \forall n\geq r,\ \forall x\in [0,\tfrac{1}{2}].
$$
Integrating on $E$ and employing \eqref{eq1.8} we reach the desired conclusion.
\qed

The $T_o$-cylinders are given by $\Delta_{[b_1]}:=\{ \frac{1}{b_1+y} : y\in [-\tfrac{1}{2},\tfrac{1}{2}]\}$ and when $r\geq 2$ by
$\Delta_{[b_1,\ldots ,b_r]} =\Delta_{[b_1]} \cap T_o^{-1} \Delta_{[b_2]} \cap \ldots \cap T_o^{-(r-1)} \Delta_{[b_r]}$.

\proof[Proof of Corollary 2]
Writing $E=E_+ \cup (-E_+)$ with $E_+ \subseteq [0,\frac{1}{2}]$, we have $\mu_o(E)=\mu (E_+)=2\mu_o (E_+)$.
Every $T_o$-cylinder $F$ is either a $T$-cylinder or the union of two $T$-cylinders, so we can assume that
$F$ is a $T$-cylinder. We have either $F\subseteq (0,\frac{1}{2}]$ or $F\subseteq [-\frac{1}{2},0)$.

Assume first $F\subseteq (0,\frac{1}{2}]$. The equality
$$
T_o^{-n}E \cap F =(T_o^{-n} E_+ \cap F) \cup (( -T_o^{-n} E_+) \cap F)=T^{-n } E_+ \cap F
$$
and Corollary 1 entail
$$
\begin{aligned}
\mu_o (T_o^{-n} E \cap F) & =\mu_o (T^{-n} E_+ \cap F)
=2\mu (T^{-n} E_+ \cap F) \\ & =2\mu(E_+)\mu(F)+O_F(q^n)=\mu(E)\mu(F) +O_F (q^n) .
\end{aligned}
$$

When $F\subseteq [-\frac{1}{2},0)$, we employ
$\mu_o (T_o^{-n}E \cap F)=\mu_o (-T_o^{-n} E \cap F) =\mu_o (T_o^{-n}E \cap (-F))$.
\qed

\section{A refinement of Rieger's bound}
Conjugating the map $T_e$ by $J:[-\frac{1}{2},\frac{1}{2}] \longrightarrow [0,1]$, $J(x):=\begin{cases} x & \text{ if $0\leq x< \frac{1}{2}$} \\
x+1 & \text{ if $-\frac{1}{2} \leq x \leq 0$} \end{cases} $ with $J^{-1}(y)=\begin{cases} y & \text{ if $0\leq y<\frac{1}{2}$} \\
y-1 & \text{ if $\frac{1}{2} \leq y\leq 1$} \end{cases}$, one gets $\widetilde{T}_e =JT_e J^{-1}:[0,1]\longrightarrow [0,1]$, which satisfies
$$
\begin{aligned}
0 < y< \frac{1}{2} \ & \Longrightarrow \ \widetilde{T}_e(y) =JT_e J^{-1} (y) = JT_e (y)=
\begin{cases} \frac{1}{y}-k & \text{ if $\frac{2}{2k+1} < y < \frac{1}{k}$, $k\geq 2$} \\
\frac{1}{y}-k+1 & \text{ if $\frac{1}{k} < y <\frac{2}{2k-1}$, $k\geq 3$}
\end{cases} . \\
\frac{1}{2} < y <1 \ & \Longrightarrow \ \widetilde{T}_e (y)=JT_e J^{-1}(y)=JT_e (y-1)=JT_e (1-y)= JT_e J^{-1}(1-y)=\widetilde{T}_e (1-y).
\end{aligned}
$$

\begin{figure}
\includegraphics[scale=0.53]{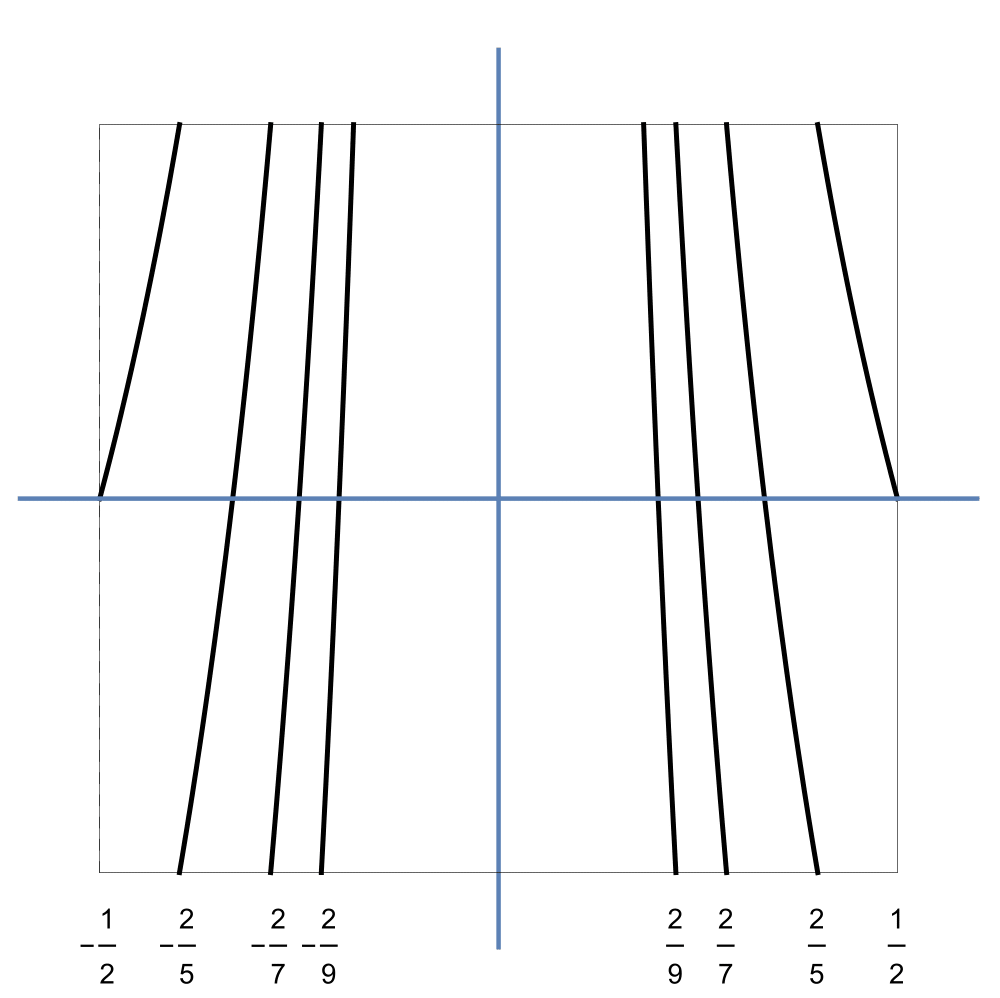}\qquad
\includegraphics[scale=0.48]{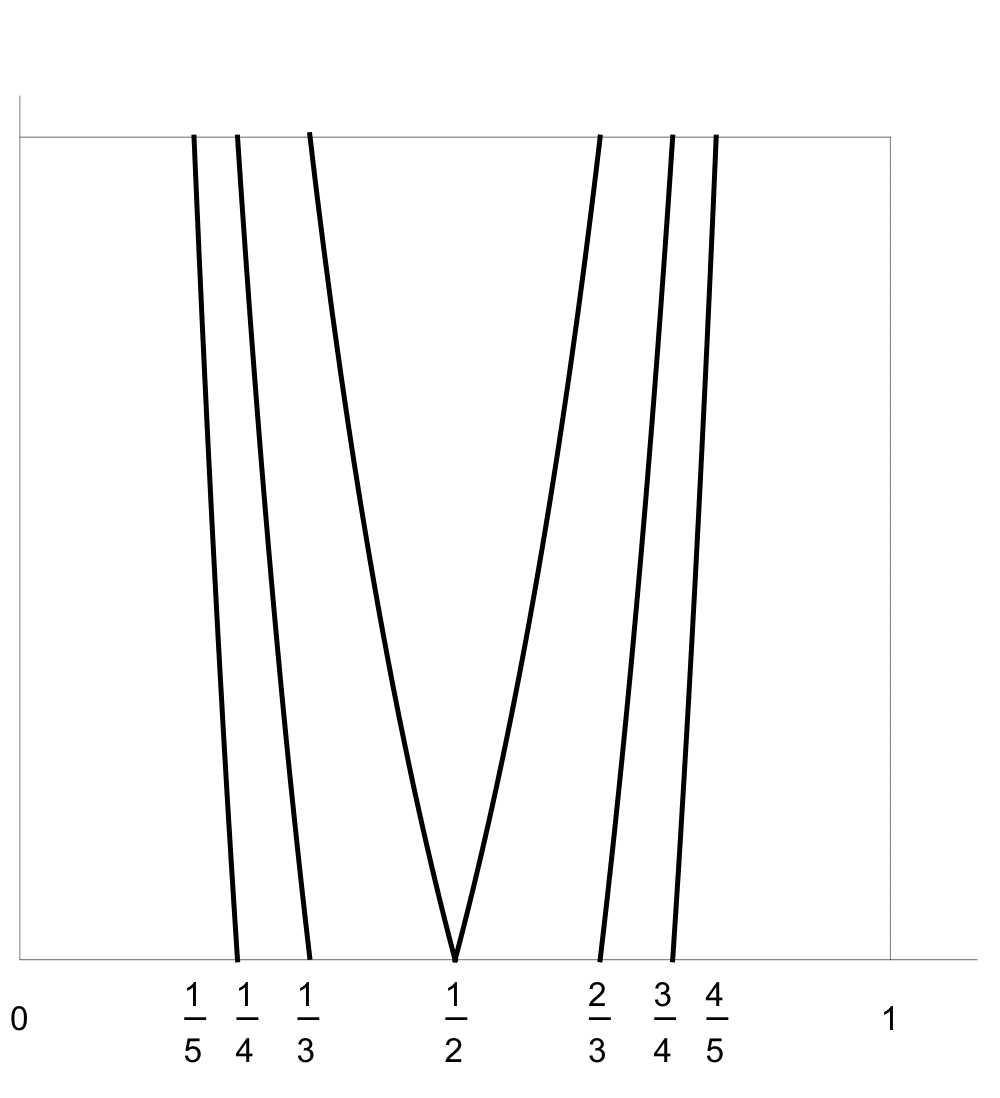}
\caption{Graphs of the maps $T_e$ and $\widetilde{T}_e$.}\label{Figure2}
\end{figure}

Observe that
$$
\widetilde{T}_e(x)=\GG (x)=\frac{1}{x}-\bigg\lfloor \frac{1}{x}\bigg\rfloor ,\quad \forall x\in ( 0,\tfrac{1}{2}].
$$
The push-forward probability measure $\widetilde{\mu}_e=J_* \mu_e =\widetilde{h}_e d\lambda$, $\widetilde{h}_e(x)=\frac{1}{(G+x)\log G}$, is $\widetilde{T}_e$-invariant.
We also have $\lambda (S) =\lambda (JS)$ for every Borel set $S\subseteq [-\frac{1}{2},\frac{1}{2}]$.

The Perron-Frobenius operator $\widetilde{P} =\widehat{(\widetilde{T}_e)}_\lambda$ associated with the transformation $\widetilde{T}_e$ and the Lebesgue measure is given by
$$
(\widetilde{P} f)(y)= \sum\limits_{x\in \widetilde{T}_e^{-1}(y)} f(x(y)) \lvert x^\prime (y)\rvert =
\sum\limits_{k=2}^\infty \Bigg( f\bigg( \frac{1}{y+k}\bigg) +f \bigg( 1-\frac{1}{y+k}\bigg)\Bigg) \frac{1}{(y+k)^2} ,\quad y\in [0,1].
$$
It satisfies $\widetilde{P}\widetilde{h}_e=\widetilde{h}_e$, emphasizing that $\widetilde{\mu}_e$ is $\widetilde{T}_e$-invariant.
Set $\widetilde{H}_e=\frac{1}{\widetilde{h}_e}$. The Perron-Frobenius operator $\widetilde{U}:=\widehat{(\widetilde{T}_e)}_{\widetilde{\mu}_e}$
associated with the $\widetilde{T}_e$-invariant measure $\widetilde{\mu}_e$, given by $\widetilde{U}=M_{\widetilde{H}_e} \widetilde{P} M_{\widetilde{h}_e}$,
is explicitly computed as
\begin{equation}\label{eq3.1}
\begin{aligned}
(\widetilde{U} f)(x) & =(x+G) \sum\limits_{k=2}^\infty
\Bigg( \frac{f(\frac{1}{x+k})}{G+\frac{1}{x+k}} +\frac{f(1-\frac{1}{x+k})}{G+1-\frac{1}{x+k}} \Bigg ) \frac{1}{(x+k)^2}\\
& =\sum\limits_{k=2}^\infty  A_k(x) f\bigg( \frac{1}{k+x}\bigg) +B_k(x) f\bigg( 1-\frac{1}{k+x}\bigg)  ,
\end{aligned}
\end{equation}
where
\begin{equation}\label{eq3.2}
\begin{aligned}
A_k(x) & =\frac{G+x}{G(k+x)(k+x+G-1)} =(G+x) \bigg( \frac{1}{k+x}-\frac{1}{k+x+G-1}\bigg)> 0,
\\ B_k(x) & =\frac{G+x}{(G+1)(k+x)(k+x+G-2)} =(G+x)\bigg( \frac{1}{k+x+G-2}-\frac{1}{k+x}\bigg) >0 .
\end{aligned}
\end{equation}
Formulas \eqref{eq3.1} and \eqref{eq3.2} define a bounded linear operator $\widetilde{U}:C[0,1] \longrightarrow C[0,1]$.
Furthermore, $\widetilde{U}1=1$ and $\widetilde{U} (C^1[0,1])\subseteq C^1[0,1]$.
We have
\begin{equation}\label{eq3.3}
\sum\limits_{k=2}^\infty A_k(x)+B_k(x)=1 \quad \text{and} \quad \sum\limits_{k=2}^\infty A_k^\prime (x)+B_k^\prime (x)=0,\quad \forall x\in [0,1].
\end{equation}
We also have
\begin{equation*}
A_k(x)+B_k(x) =\frac{G+x}{(k+x+G-1)(k+x+G-2)},\quad
 \sum\limits_{k=5}^\infty A_k(x)+B_k(x)  =\frac{G+x}{G+3+x}.
 \end{equation*}

We write
$(\widetilde{U} f)^\prime =S_I f+S_{II} f$, with
$$
\begin{aligned}
(S_I f)(x) & =-\sum\limits_{k=2}^\infty \frac{A_k(x)}{(k+x)^2} \, f^\prime \bigg( \frac{1}{k+x}\bigg)
+\sum\limits_{k=2}^\infty \frac{B_k (x)}{(k+x)^2} \, f^\prime \bigg( 1-\frac{1}{k+x}\bigg), \\ (S_{II}f)(x) &
=\sum\limits_{k=2}^\infty A_k^\prime (x) f\bigg( \frac{1}{k+x}\bigg) +B_k^\prime (x) f\bigg( 1- \frac{1}{k+x}\bigg) ,\quad \forall f\in C^1[0,1],\ x\in [0,1].
\end{aligned}
$$

\begin{Lem2*}
$\| S_I f \|_\infty \leq  0.1346 \| f^\prime\|_\infty$, $\forall f\in C^1[0,1]$.
\end{Lem2*}

\proof[Proof.]
We have $\lvert (S_I f)(x)\rvert \leq \Phi (x)\| f^\prime \|_\infty$, where $\Phi=\Phi_2+\Phi_3+\Phi_4+\Phi_5$, with
$$\begin{aligned}
\Phi_2 (x) & =\frac{A_2(x)+B_2(x)}{(2+x)^2} = \frac{1}{(2+x)^2 (G+1+x)} , \\
\Phi_3 (x) & =\frac{A_3 (x)+B_3(x)}{(3+x)^2} = \frac{G+x}{(3+x)^2(G+1+x)(G+2+x)}  ,\\
\Phi_4 (x) & =\frac{A_4 (x)+B_4(x)}{(4+x)^2} = \frac{G+x}{(4+x)^2(G+2+x)(G+3+x)}   ,\\
\Phi_5(x) & = \frac{1}{(5+x)^2} \sum\limits_{k=5}^\infty A_k(x)+B_k(x) =\frac{G+x}{(5+x)^2 (G+3+x)}.
\end{aligned}
$$
The function $\Phi$ is decreasing on $[0,1]$ with
$\| \Phi\|_\infty \leq \Phi(0) <  0.1346.
$
\qed

\begin{Lem3*}
$\| S_{II} f\|_\infty \leq 0.092 \| f^\prime\|_\infty$, $\forall f\in C^1[0,1]$.
\end{Lem3*}

\proof[Proof.]
Compute
$$
\begin{aligned}
A_k^\prime (x) & = \bigg( \frac{1}{k+x}-\frac{1}{k+x+G-1}\bigg) \Bigg( 1-(G+x)\bigg( \frac{1}{k+x}+\frac{1}{k+x+G-1}\bigg) \Bigg),  \\
B_k^\prime (x) & = \bigg( \frac{1}{k+x+G-2}-\frac{1}{k+x}\bigg) \Bigg( 1-(G+x)\bigg( \frac{1}{k+x+G-2}+\frac{1}{k+x}\bigg) \Bigg) , \\
A_k^\prime (x)+B_k^\prime (x) & = \frac{1}{k+x+G-2}-\frac{1}{k+x+G-1} \\
& \hspace{1cm} -(G+x) \bigg( \frac{1}{(k+x+G-2)^2} -\frac{1}{(k+x+G-1)^2} \bigg) .
\end{aligned}
$$
Using the second identity in \eqref{eq3.3} we can write
$$
(S_{II} f)(x)=\sum\limits_{k=2}^\infty A_k^\prime (x) \Bigg( f\bigg( \frac{1}{k+x}\bigg) -f\bigg( \frac{1}{2}\bigg)\Bigg)
+B_k^\prime (x) \Bigg( f\bigg( 1-\frac{1}{k+x}\bigg) -f\bigg( \frac{1}{2}\bigg)\Bigg) .
$$
In conjunction with the Mean Value Theorem and $\lvert \frac{1}{2+x} -\frac{1}{2}\rvert\leq \frac{1}{6}$,
$\lvert \frac{1}{3+x} -\frac{1}{2}\rvert \leq \frac{1}{4}$,
$\lvert \frac{1}{4+x} -\frac{1}{2}\rvert\leq \frac{3}{10}$,
$\lvert \frac{1}{k+x} -\frac{1}{2}\rvert\leq \frac{1}{2}$, $k\geq 5$,
this yields
$$
\lvert (S_{II} f)(x)\rvert \leq \Psi (x) \| f^\prime\|_\infty,
$$
with $\Psi=\Psi_2+\Psi_3+\Psi_4+\Psi_5$, where
$$
\Psi_2 =\frac{\lvert A_2^\prime\rvert+\lvert B_2^\prime\rvert}{6},\quad
\Psi_3 =\frac{\lvert A_3^\prime\rvert+\lvert B_3^\prime \rvert}{4}, \quad
\Psi_4 =  \frac{3(\lvert A_4^\prime\rvert+\lvert B_4^\prime \rvert)}{10}, \quad
\Psi_5 = \frac{1}{2} \sum\limits_{k=5}^\infty
\lvert A_k^\prime\rvert+\lvert B_k^\prime\rvert .
$$

When $k\geq 5$ we have $A_k^\prime >0$ and $B_k^\prime >0$ on $[0,1]$. The above expression for $A_k^\prime+B_k^\prime$ allows us to compute
$$
\begin{aligned}
\Psi_5(x) & =\frac{1}{2}\sum\limits_{k=5}^\infty
A_k^\prime(x)+ B_k^\prime(x) =\frac{1}{2(G+3+x)}-\frac{G+x}{2(G+3+x)^2} \\ & =\frac{3}{2(G+3+x)^2}
\leq \frac{3}{2(G+3)^2} < 0.0704.
\end{aligned}
$$
On the other hand we have $A_2^\prime <0$ and $B_2^\prime <0$ on $[0,1]$, leading to
$$
\begin{aligned}
\Psi_2 (x) & =-\frac{A_2^\prime(x)+B_2^\prime (x)}{6}
=\frac{1}{6}\Bigg( \frac{1}{G+1+x}-\frac{1}{G+x} +(G+x)\bigg( \frac{1}{(G+x)^2}-\frac{1}{(G+1+x)^2}\bigg) \Bigg) \\
& =\frac{1}{6(G+1+x)^2} \leq \frac{1}{6(G+1)^2}=\frac{(2-G)^2}{6} <0.0244.
\end{aligned}
$$
Numerically, we see that $\Psi_3(x) \leq \frac{1}{4} (\lvert A_3^\prime (1)\rvert+\lvert B_3^\prime (1)\rvert) <0.0019$
and $\Psi_4(x)\leq \frac{3}{10} (\lvert A_4^\prime (0)\rvert+\lvert B_4^\prime (0)\rvert)<0.0025$. Thus
$\| S_{II} f\|_\infty \leq 0.0992 \| f^\prime \|_\infty$.
\qed

\begin{Cor3*}
$\| (\widetilde{U} f)^\prime \|_\infty \leq 0.234 \| f^\prime \|_\infty$, $\forall f\in C^1[0,1]$.
\end{Cor3*}

The proof of the following asymptotic formula follows ad litteram the proof of Theorem 1 (i).

\begin{Prop2*}
With $q=0.234$, for every Borel set $E\subseteq [0,1]$,
$$
\lambda (\widetilde{T}_e^{-n} E)=\widetilde{\mu}_e (E) +O(\widetilde{\mu}_e (E)q^n).
$$
\end{Prop2*}

Let $E\subseteq [-\frac{1}{2},\frac{1}{2}]$ be a Borel set and $\widetilde{E}:=JE$.
The equality $JT_e =\widetilde{T}_e J$ yields $JT_e^{-n} E=\widetilde{T}_e^{-n} JE$.
Theorem 1 (iii) now follows from $\lambda (JS)=\lambda (S)$ for every Borel set $S\subseteq [-\frac{1}{2},\frac{1}{2}]$, Proposition 2, $\widetilde{\mu}_e(\widetilde{E})=\mu_e(E)$, and
\begin{equation*}
\begin{aligned}
\lambda (T_e^{-n} E) & =\lambda (J^{-1}\widetilde{T}_e^{-n} J E) =\lambda (\widetilde{T}_e^{-n} \widetilde{E}) \\ &
= \widetilde{\mu}_e (\widetilde{E}) +O(\widetilde{\mu}_e (\widetilde{E}) q^n) =\mu_e (E)+O(\mu_e (E) q^n).
\end{aligned}
\end{equation*}

The $T_e$-cylinders are given (up to null sets) by
$\Delta^e_{[\![(a_1,\pm 1)]\!]} =\pm  (\frac{2}{2a_1+1},\frac{2}{2a_1-1})$, $a_1 \geq 3$,
$\Delta^e_{[\![(2,\pm 1)]\!]} =\pm (\frac{2}{5},\frac{1}{2})$, and when $r\geq 2$ by
$\Delta^e_{[\![(a_1,e_1),\ldots,(a_r,e_r)]\!]} =\Delta^e_{[\![(a_1,e_1)]\!]} \cap T_e^{-1} \Delta^e_{[\![(a_2,e_2)]\!]}
\cap \ldots \cap  T_e^{-(r-1)} \Delta^e_{[\![(a_r,e_r)]\!]}$.

The $\widetilde{T}_e$-cylinders are given by $\widetilde{\Delta}^e_{[\![(a_1,+1)]\!]}=(\frac{1}{a_1+1},\frac{1}{a_1})$,
$\widetilde{\Delta}^e_{[\![(a_1,-1)]\!]} =1-\widetilde{\Delta}^e_{[\![(a_1,+1)]\!]}$, $a_1 \geq 2$,
and when $r\geq 2$ by
$\widetilde{\Delta}^e_{[\![(a_1,e_1),\ldots,(a_r,e_r)]\!]} =\widetilde{\Delta}^e_{[\![(a_1,e_1)]\!]} \cap \widetilde{T}_e^{-1} \widetilde{\Delta}^e_{[\![(a_2,e_2)]\!]}
\cap \ldots \cap  \widetilde{T}_e^{-(r-1)} \widetilde{\Delta}^e_{[\![(a_r,e_r)]\!]}$.

Note that $J$ does not map $T_e$-cylinders into $\widetilde{T}_e$-cylinders.

Formula \eqref{eq3.1} shows in particular that $\widetilde{U}$ acts as
$$
(\widetilde{U} f)(x)=\sum\limits_{\substack{k\geq 2 \\ e=\pm 1}} \widetilde{P}_{(k,e)} (x) f (\widetilde{w}_{(k,e)} (x)),
$$
where $\widetilde{w}_{(k,+1)}(x)=\frac{1}{k+x}$, $\widetilde{w}_{(k,-1)}(x)=1-\frac{1}{k+x}$ map the interval $[0,1]$
onto the rank one cylinders $\widetilde{\Delta}^e_{[\![(k,\pm 1)]\!]}$ respectively.
As a result, the argument in the proof of Corollary 1 applies, entailing

\begin{equation}\label{eq3.4}
\widetilde{\mu}_e( \widetilde{T}_e^{-n} \widetilde{E} \cap \widetilde{F})=\widetilde{\mu}_e(\widetilde{E}) \widetilde{\mu}_e(\widetilde{F})+O_{\widetilde{F}} (q^n),
\end{equation}
for any Borel set $\widetilde{E} \subseteq [0,1]$ and any $\widetilde{T}_e$-cylinder $\widetilde{F}$, with $q=0.234$.

\begin{Cor4*}
With $q=0.234$, for any Borel set $E\subseteq [-\frac{1}{2},\frac{1}{2}]$ and $F=J^{-1} \widetilde{F} \subseteq [-\frac{1}{2},\frac{1}{2}]$,
$\widetilde{F}$ $\widetilde{T}_e$-cylinder,
$$
\mu_e(T_e^{-n} E \cap F)=\mu_e(E) \mu_e(F)+O_F (q^n).
$$
\end{Cor4*}

\proof
Employing $\mu_e (S)=\widetilde{\mu}_e (JS)$ and $JT_e^{-n} S=\widetilde{T}_e^{-n} JS$ for every Borel set $S\subseteq [-\frac{1}{2},\frac{1}{2}]$, and equation \eqref{eq3.4} with $\widetilde{E}:=JE$,  we find
\begin{equation*}
\begin{aligned}
\mu_e (T_e^{-n}E\cap F) & =\widetilde{\mu}_e (J T_e^{-n} E \cap J F) = \widetilde{\mu}_e (\widetilde{T}_e^{-n} \widetilde{E} \cap \widetilde{F})
\\
& = \widetilde{\mu}_e (\widetilde{E})\widetilde{\mu}_e (\widetilde{F}) +O_F (q^n) =\mu_e (E) \mu_e (F) +O_F(q^n) ,
\end{aligned}
\end{equation*}
which concludes the proof.
\qed

\bigskip

{\bf Acknowledgments.} This research was supported by NSF award DMS-1449269 and University of Illinois Research Board Award RB-22069.

We are grateful to Joseph Vandehey for constructive comments and to the referees for careful reading and valuable suggestions.

\end{document}